\documentclass[10pt,reqno,draft]{amsart}
\usepackage[a4paper,centering,hscale=0.7,vscale=0.75]{geometry}
\usepackage[T1]{fontenc}
\usepackage{mathtools}

\usepackage{amssymb}
\usepackage{amsmath}
\usepackage{amsthm}
\usepackage{mathrsfs}
\usepackage{bbm}

\newtheorem{thm}{Theorem}[section]

\newtheorem{lem}[thm]{Lemma}
\newtheorem{prop}[thm]{Proposition}
\newtheorem{defi}[thm]{Definition}

\newcommand{\ind}[1]{\mathbbm{1}_{#1}}

\newcommand{\dom}{\mathsf{D}}
\newcommand{\cF}{\mathscr{F}}
\newcommand{\cL}{\mathscr{L}}
\newcommand{\enne}{\mathbb{N}}
\renewcommand{\P}{\mathbb{P}}
\newcommand{\erre}{\mathbb{R}}
\newcommand{\cR}{\mathscr{R}}

\DeclarePairedDelimiter{\abs}{\lvert}{\rvert}
\DeclarePairedDelimiter{\norm}{\lVert}{\rVert}
\DeclarePairedDelimiterX\ip[2]{\langle}{\rangle}{#1,#2}

\DeclarePairedDelimiterX\cc[2]{[\![}{]\!]}{#1,#2}
\DeclarePairedDelimiterX\co[2]{[\![}{[\![}{#1,#2}

\def\OO{\mathcal{O}}

\begin{document}
\title[Positivity of mild solutions]{On the positivity of local mild
  solutions to stochastic evolution equations}

\author{Carlo Marinelli}
\address[Carlo Marinelli]{Department of Mathematics, University College London, 
Gower Street, London WC1E 6BT, United Kingdom}
\urladdr{https://www.homepages.ucl.ac.uk/$\sim$ucahcm0/}

\author{Luca Scarpa}
\address[Luca Scarpa]{Faculty of Mathematics, University of Vienna, 
Oskar-Morgenstern-Platz 1, 1090 Wien, Austria.}
\email{luca.scarpa@univie.ac.at}
\urladdr{http://www.mat.univie.ac.at/$\sim$scarpa}

\date{December 12, 2019}



\begin{abstract}
  We provide sufficient conditions on the coefficients of a stochastic
  evolution equation on a Hilbert space of functions driven by a
  cylindrical Wiener process ensuring that its mild solution is
  positive if the initial datum is positive. As an application, we
  discuss the positivity of forward rates in the Heath-Jarrow-Morton
  model via Musiela's stochastic PDE.
\end{abstract}

\maketitle


\section{Introduction}
\label{sec:intro}
Let us consider a stochastic evolution equation of the type
\begin{equation}
  \label{eq:0}
  du + Au\,dt = F(u)\,dt + B(u)\,dW, \qquad u(0)=u_0,
\end{equation}
where $A$ is a linear maximal monotone operator on a Hilbert space of
functions $H$, the coefficients $F$ and $B$ satisfy suitable
integrability assumptions, and $W$ is a cylindrical Wiener process.
Precise assumptions on the data of the Cauchy problem \eqref{eq:0} are
given in \S\ref{sec:main} below. Our goal is to establish a maximum
principle for (local) mild solutions to \eqref{eq:0}, i.e. to provide
sufficient conditions on the operator $A$ and on the coefficients $F$
and $B$ such that positivity of the initial datum $u_0$ implies
positivity of the solution $u$ (see Theorem~\ref{thm:pos} below).

A simpler problem was studied in \cite{cm:pos1}, where coefficients $F$
and $B$ are assumed to be Lipschitz continuous. Here we simply assume
that $F$ and $B$ satisfy rather minimal integrability conditions and
that a local mild solution exists. On the other hand, in \cite{cm:pos1}
the linear operator $A$ need only generate a positivity preserving
semigroup, while here we require that $A$ generates a
sub-Markovian semigroup.

We refer to \cite{cm:pos1} for a discussion about the relation of other
positivity results for solutions to stochastic partial differential
equations with ours. It is however probably worth pointing out that
most existing results seem to deal with equations in the variational
setting (see, e.g., \cite{Kry:MP-SPDE,Kry:shortIto,Pard}).

As an application, we provide an alternative, more direct proof of the
positivity of forward rates in the Heath-Jarrow-Morton \cite{HJM}
framework with respect to the one in \cite{cm:pos1}. This is obtained, as
is now classical, viewing forward curves as solutions to the so-called
Musiela stochastic PDE (see, e.g., \cite{filipo,cm:MF10}).


\section{Assumptions and main result}
\label{sec:main}
Let $(\Omega,\cF,\P)$ be a probability space endowed with a complete
right-continuous filtration $(\cF_t)_{t\in[0,T]}$, with $T>0$ a fixed
final time, on which all random elements will be defined. Identities
and inequalities between random variables are meant to hold
$\P$-almost surely, and two stochastic processes are declared equal,
unless otherwise stated, if they are indistinguishable.  The
$\sigma$-algebra of progressively measurable subsets of
$\Omega\times[0,T]$ will be denoted by $\cR$.
We shall denote a cylindrical Wiener process on a separable Hilbert
space $U$ by $W$.
Standard notation and terminology of stochastic calculus for
semimartingales will be used throughout (see, e.g., \cite{Met}). In
particular, given an adapted process $X$ and a stopping time $\tau$,
$X^\tau$ will denote the process $X$ stopped at $\tau$.  Similarly, if
$X$ is also c\`adl\`ag, $X^{\tau-}$ stands for the process $X$
pre-stopped at $\tau$.

For any separable Hilbert spaces $E_1$ and $E_2$, we will use the
symbols $\cL(E_1,E_2)$ $\cL^2(E_1,E_2)$ for the space of linear
continuous and Hilbert-Schmidt operators from $E_1$ to $E_2$,
respectively. The space of continuous bilinear maps from
$E_1 \times E_1$ to $E_2$ will be denoted by $\cL_2(E_1;E_2)$.
The $n$-th order Fr\'echet and G\^ateaux derivatives of a function
$\Phi: E_1 \to E_2$ at a point $x \in E_1$ are denoted by $D^n\Phi(x)$
and $D^n_{\mathcal{G}}\Phi(x)$, respectively, omitting the superscript
if $n=1$, as usual.

\medskip

We shall work under the following standing assumptions.
\smallskip\par\noindent
\textbf{(A1)} There exists an open set $\mathcal{O}$ in $\erre^d$,
$d \geq 1$, and a Borel measure $\mu$ such that
$H=L^2(\mathcal{O},\mu)$.
\smallskip\par\noindent
The norm and scalar product on $H$ will be denoted by $\norm{\cdot}$
and $\ip{\cdot}{\cdot}$, respectively.
\smallskip\par\noindent
\textbf{(A2)} $A$ is a linear maximal monotone operator on $H$ such
that its resolvent is sub-Markovian and is a contraction with respect
to the $L^1(\OO,\mu)$-norm.
\smallskip\par\noindent
Recall that the resolvent of $A$, i.e. the family of linear continuous
operators on $H$ defined by
\[
J_\lambda := (I+\lambda A)^{-1}, \qquad \lambda>0,
\]
is said to be sub-Markovian if, for every $\lambda>0$ and every
$\phi \in H$ such that $0 \leq \phi \leq 1$ a.e. in $\mathcal{O}$, one has
$0 \leq J_\lambda\phi \leq 1$ a.e. in $\mathcal{O}$.
\smallskip\par\noindent
\textbf{(A3)} $F:\Omega \times[0,T] \times H \to H$ and
$B:\Omega \times[0,T] \times H \to \cL^2(U,H)$ are
$\cR \otimes \mathscr{B}(H)$-measurable, and there 
exists a constant $C>0$ such that 
\[ 
  -\ip{F(\omega,t,h)}{h_-}
  + \frac12\norm[\big]{1_{\{h<0\}}B(\omega,t,h)}_{\cL^2(U,H)}^2
  \leq C \norm{h_-}^2_{L^2(\OO)}
  \qquad \forall (\omega,t,h) \in \Omega \times [0,T] \times H.
\]
In particular, note that choosing $h=0$ yields $F(\cdot, 0)=0$ and
$B(\cdot, 0)=0$.
\smallskip\par\noindent
\textbf{(A4)} $u_0\in L^0(\Omega,\cF_0; H)$

\medskip

\begin{defi}
  A local mild solution to the Cauchy problem \eqref{eq:0} is a pair
  $(u,\tau)$, where $\tau$ is a stopping time with $\tau \leq T$, and
  $u:[\![0,\tau[\![ \to H$ is a measurable adapted process with
  continuous trajectories such that, for any stopping time
  $\sigma<\tau$, one has
  \begin{itemize}
  \item[(i)] $S(t-\cdot) F(u) \ind{\cc{0}{\sigma}} \in L^0(\Omega;L^1(0,t;H))$
    for all $t \in [0,T]$;
  \item[(ii)] $S(t-\cdot) B(u) \ind{\cc{0}{\sigma}}
    \in L^0(\Omega;L^2(0,t;\cL^2(U,H)))$ for all $t \in [0,T]$,
  \end{itemize}
  and
  \[
    u = S(\cdot)u_0 + \int_0^\cdot S(\cdot-s)F(s,u(s))\,ds +
    \int_0^\cdot S(\cdot-s)B(s,u(s))\,dW(s).
  \]
\end{defi}
The last identity is to be understood in the sense of
indistinguishability of processes defined on the stochastic interval
$\co{0}{\tau}$. Here the stochastic convolution is defined on
$\cc{0}{\sigma}$, for every stopping time $\sigma<\tau$, as
\[
  \biggl( \int_0^t S(t-s) B(s,u(s)) \ind{\cc{0}{\sigma}}(s)\,dW(s)
  \biggr)_{t\in[0,\sigma]}.
\]

The main result is the following.
\begin{thm}
  \label{thm:pos}
  Let $(u,\tau)$ be a local mild solution to the Cauchy problem
  \eqref{eq:0} such that, for every stopping time $\sigma<\tau$, one
  has
  \begin{itemize}
  \item[(i)] $F(u)\ind{\cc{0}{\sigma}} \in L^0(\Omega;L^1(0,T;H))$;
  \item[(ii)] $B(u)\ind{\cc{0}{\sigma}} \in L^0(\Omega;L^2(0,T;\cL^2(U,H)))$.
  \end{itemize}
  If $u_0 \geq 0$ a.e. in $\OO$, then $u^{\tau-}(t) \geq 0$ a.e.~in
  $\OO$ for all $t \in [0,T]$.
\end{thm}


\section{Auxiliary results}
\label{sec:aux}
The arguments used in the proof of Theorem~\ref{thm:pos} (see
\S\ref{sec:proof} below) rely on the following results, that we recall
here for the reader's convenience.
The first is a continuous dependence result for mild solutions to
stochastic evolution equations in the form \eqref{eq:0} with respect
to the coefficients and the initial datum. This is a consequence of a
more general statement proved in \cite[Corollary~3.4]{KvN2}. Let
\begin{align*}
  (u_{0n})_n &\subset L^0(\Omega,\cF_0;H),\\
  (f_n)_n, f &\subset L^0(\Omega;L^1(0,T;H)),\\
  (G_n)_n, G &\subset L^0(\Omega;L^2(0,T;\cL^2(U,H)))
\end{align*}
be such that the $H$-valued processes $f_n$, $f$, $G_nv$, and $Gv$ are
strongly measurable and adapted for all $v \in U$ and $n \in \enne$.
Then the Cauchy problems
\[
  du_n + Au_n\,dt = f_n\,dt + G_n\,dW, \qquad u_n(0)=u_{0n},
\]
and
\[
  du + Au\,dt = f\,dt + G\,dW, \qquad u(0)=u_0,
\]
admit unique mild solutions $u_n$ and $u$, respectively.
\begin{prop}
  \label{prop:micia}
  Assume that
  \begin{alignat*}{3}
  u_{0n} &\longrightarrow u_0  &&\quad \text{in } L^0(\Omega; H),\\
  f_n &\longrightarrow f &&\quad \text{in } L^0(\Omega;L^1(0,T;H)),\\
  G_n &\longrightarrow G &&\quad \text{in } L^0(\Omega; L^2(0,T;\cL^2(U,H))).
  \end{alignat*}
  Then $u_n \to  u$ in $L^0(\Omega;C([0,T];H))$.
\end{prop}

The second result we shall need is a generalized It\^o formula, the
proof of which can be found in \cite{cm:pos1}.
\begin{prop}
  \label{prop:Ito}
  Let $G \colon H \to \erre$ be continuously Fr\'echet differentiable
  and $DG$ be G\^ateaux differentiable, with
  $D^2_{\mathcal{G}}G \colon H \to \cL_2(H;\erre)$ such that
  $(\varphi,\zeta_1,\zeta_2) \mapsto
  D^2_{\mathcal{G}}G(\varphi)[\zeta_1,\zeta_2]$ is continuous, and
  assume that $G$, $DG$, and $D^2_{\mathcal{G}}G$ are polynomially
  bounded.  Moreover, let $f \in L^0(\Omega;L^1(0,T;H))$ and
  $\Phi \in L^0(\Omega;L^2(0,T;\cL^2(U,H)))$ be measurable and adapted
  processes, and $v_0 \in L^0(\Omega,\cF_0;H)$. Setting
  \[
  v := v_0 + \int_0^\cdot f(s)\,ds + \int_0^\cdot \Phi(s)\,dW(s),
  \]
  one has  
  \begin{align*}
    G(v) &= G(v_0) + \int_0^\cdot \Bigl( DG(v)f %
    + \frac12 \operatorname{Tr}\bigl(
    \Phi^* D^2_{\mathcal{G}}G(v)\Phi \bigr)\Bigr)(s)\,ds\\
    &\quad + \int_0^\cdot DG(v(s)) \Phi(s)\,dW(s).
  \end{align*}
\end{prop}

Finally, we recall an inequality for maximal monotone linear operators
with sub-Markovian resolvent, due to Br\'ezis and Strauss (see
\cite[Lemma~2]{BreStr}).\footnote{For a related inequality cf. also
  \cite[Lemma~5.1]{RW:nonmon}.}
\begin{lem}
  \label{lem:brez-str}
  Let $\beta \colon \erre \to 2^\erre$ be a maximal monotone graph
  with $0 \in \beta(0)$.  Let $\varphi \in L^p(\OO)$ with
  $A\varphi \in L^p(\OO)$, and $z \in L^q(\OO)$ with
  $z\in\beta(\varphi)$ a.e.~in $\OO$, where $p,q \in [1,+\infty]$ and
  $1/p+1/q=1$.  Then
  \[
  \int_\OO (A\varphi) z \geq 0.
  \]
\end{lem}
We include a sketch of proof for the reader's convenience, assuming
for simplicity that $\beta: \erre \to \erre$ is continuous and bounded. Let
$j:\erre \to \erre_+$ a (differentiable, convex) primitive of $\beta$
and
\[
  A_\lambda := \frac{1}{\lambda}\bigl( I - (I+\lambda A)^{-1}\bigr)
  = \frac{1}{\lambda}(I-J_\lambda), \qquad \lambda>0,
\]
the Yosida approximation of $A$. It is well known that $A_\lambda$ is
a linear maximal monotone bounded operator on $H$ and that, for every
$v \in \dom(A)$, $A_\lambda v \to Av$ as $\lambda \to 0$. Let
$v \in \dom(A)$. The convexity of $j$ implies, for every $\lambda>0$,
\begin{align*}
  \ip[\big]{A_\lambda v}{\beta(v)}_{L^2}
  &= \frac{1}{\lambda} \ip[\big]{v - J_\lambda v}{j'(v)}_{L^2}\\
  &\geq \frac{1}{\lambda} \biggl(
    \int_\OO j(v) - \int_\OO j(J_\lambda v) \biggr)
    = \frac{1}{\lambda} \bigl( \norm{j(v)}_{L^1}
    - \norm{j(J_\lambda v)}_{L^1} \bigr).
\end{align*}
Since $J_\lambda$ is sub-Markovian and $j$ is convex, the generalized
Jensen inequality for positive operators (see \cite{Haa07}) and the
contractivity of $J_\lambda$ in $L^1$ imply that
\[
  \norm[\big]{j(J_\lambda v)}_{L^1} \leq \norm[\big]{J_\lambda j(v)}_{L^1}
  \leq \norm[\big]{j(v)}_{L^1},
\]
i.e. that
\[
\ip[\big]{A_\lambda v}{\beta(v)}_{L^2} \geq 0
\]
for every $\lambda \to 0$. Passing to the limit as $\lambda \to 0$
yields $\ip{Av}{\beta(v)}_{L^2} \geq 0$.


\section{Proof of Theorem~\ref{thm:pos}}
\label{sec:proof}
The proof is divided into two parts. First we show that a local mild
solution $u$ to \eqref{eq:0} can be approximated by strong solutions
to regularized equations. As a second step, we show that such
approximating processes are positive, thanks to a suitable version of
It\^o's formula.

\subsection{Approximation of the solution}
Let $(u,\tau)$ be a local mild solution to \eqref{eq:0}.  Let $\sigma$
be a stopping time with $\sigma<\tau$, so that
$u:[\![0,\sigma]\!]\to H$ is well defined, and set
\begin{align*}
  \bar{u} &:= u^\sigma \in L^0(\Omega;C([0,T];H)),\\
  \bar{F} &:= F(\cdot,u) \ind{\cc{0}{\sigma}}
            \in L^0(\Omega;L^1(0,T;H)),\\
  \bar{B} &:= B(\cdot,u) \ind{\cc{0}{\sigma}}
            \in L^0(\Omega;L^2(0,T;\cL^2(U,H))).
\end{align*}
Note that, by assumption (A3), $F(\cdot,0)=0$ and $B(\cdot,0)=0$, hence
\begin{align*}
  \bar{F} = F(\cdot,u) \ind{\cc{0}{\sigma}}
  = F(\cdot,u\ind{\cc{0}{\sigma}}), \\
  \bar{B} = B(\cdot,u)\ind{\cc{0}{\sigma}}
  =B(\cdot,u\ind{\cc{0}{\sigma}}).
\end{align*}
In particular, one has
\begin{equation}
  \label{mild_bar}
  \bar{u}(t) := S(t)u_0
  + \int_0^t S(t-s)\bar F(s)\,ds 
  + \int_0^t S(t-s)\bar B(s)\,dW(s)
\end{equation}
for all $t \in [0,T]$ $\P$-a.s., or, equivalently, $\bar{u}$ is the
unique global mild solution to the Cauchy problem
\[
  d\bar{u} + A\bar{u}\,dt =\bar{F}\,dt + \bar{B}\,dW,
  \qquad \bar{u}(0)=u_0.
\]
Recalling that $J_\lambda \in \cL(H,\dom(A))$ for all $\lambda>0$, one has
\begin{align*}
  \bar{F}_\lambda
  &:= J_\lambda F(\cdot,u) \ind{\cc{0}{\sigma}}
    = J_\lambda \bar{F} \in L^0(\Omega;L^1(0,T;\dom(A))),\\
  \bar{B}_\lambda
  &:=J_\lambda B(\cdot,u) \ind{\cc{0}{\sigma}}
    = J_\lambda \bar{B} \in L^0(\Omega;L^2(0,T;\cL^2(U,\dom(A)))),\\
  u_{0\lambda}
  &:= J_\lambda u_0 \in L^0(\Omega,\cF_0;\dom(A)),
\end{align*}
where the second assertion is an immediate consequence of the ideal
property of Hilbert-Schmidt operators.
The process $u_\lambda:\Omega\times[0,T]\to H$ defined as
\begin{equation}
  \label{mild_bar_lam}
  u_\lambda(t) := S(t)u_{0\lambda }
  + \int_0^t S(t-s)\bar F_\lambda(s)\,ds 
  + \int_0^t S(t-s)\bar B_\lambda(s)\,dW(s), 
  \qquad t\in[0,T],
\end{equation}
therefore belongs to $L^0(\Omega;C([0,T];\dom(A)))$ and is the unique
global strong solution to the Cauchy problem
\[
  du_\lambda + Au_\lambda\,dt = \bar{F}_\lambda\,dt
  + \bar{B}_\lambda\,dW, \qquad u_\lambda(0)=u_{0\lambda},
\]
i.e.
\begin{equation}
  \label{strong}
  u_\lambda + \int_0^\cdot Au_\lambda(s)\,ds =
  u_{0\lambda} + \int_0^\cdot \bar{F}_\lambda(s)\,ds 
  + \int_0^\cdot \bar{B}_\lambda(s)\,dW(s)
\end{equation}
in the sense of indistinguishable $H$-valued processes.
Furthermore, since $J_\lambda$ is contractive and converges to the
identity in the strong operator topology of $\cL(H,H)$ as
$\lambda \to 0$, i.e. $J_\lambda h \to h$ for every $h \in H$, 
one has
\begin{align*}
  u_{0\lambda} \longrightarrow u_0
  &\quad \text{in } L^0(\Omega; H),\\
  \bar{F}_\lambda \longrightarrow \bar{F}
  &\quad \text{in } L^0(\Omega; L^2(0,T; H)),\\
  \bar{B}_\lambda \longrightarrow \bar{B}
  &\quad \text{in } L^0(\Omega;L^2(0,T;\cL^2(U,H))),
\end{align*}
where the second convergence follows immediately by the dominated
convergence theorem, and the third one by a continuity property of
Hilbert-Schmidt operators (see, e.g., \cite[Theorem~9.1.14]{HvNVW2}).
Finally, thanks to Proposition~\ref{prop:micia}, we deduce that 
\begin{equation}
  \label{conv}
  u_\lambda \longrightarrow \bar{u} \quad\text{in } L^0(\Omega;C([0,T];H)).
\end{equation}

\subsection{Positivity}
Let us introduce the functional 
\begin{align*}
  G \colon H
  &\longrightarrow \erre_+,\\
  G \colon \varphi
  &\longmapsto \frac12 \int_\OO  \abs{\varphi_-}^2,
\end{align*}
as well as the family, indexed by $n \in \enne$, of regularized
functionals
\begin{align*}
  G_n \colon H &\longrightarrow \erre_+,\\
  G_n \colon \varphi &\longmapsto \frac12 \int_\OO  g_n(\varphi),
\end{align*}
where $g_n:\erre \to \erre_+$ is convex, twice continuously
differentiable, identically equal to zero on $\erre_+$, strictly
positive and decreasing on $\erre_-$, such that $(g_n'')$ is uniformly
bounded, and $g'_n(r) \to -r^-$ as $n \to \infty$ for all
$r \in \erre$. The existence of such an approximating sequence is well
known (see, e.g., \cite[\S3]{scar-stef-order} for details).
One can verify (see, e.g., \cite{cm:pos1}) that, for every $n \in \enne$,
$G_n$ is everywhere continuously Fr\'echet differentiable with
derivative
\begin{align*}
  DG_n \colon H &\longrightarrow \cL(H,\erre) \simeq H,\\
  DG_n \colon \varphi &\longmapsto g_n'(\varphi),
\end{align*}
and that $DG_n \colon H \to H$ is G\^ateaux differentiable with
G\^ateaux derivative given by
\begin{align*}
  D^2_{\mathcal G}G_n \colon H &\longrightarrow \cL(H,H) \simeq \cL_2(H;\erre),\\
  D^2_{\mathcal G}G_n \colon \varphi &\longmapsto \Bigl[
  (\zeta_1,\zeta_2) \mapsto \int_{\OO} g_n''(\varphi)\zeta_1\zeta_2 \Bigr].
\end{align*}
Furthermore, the map
$(\varphi,\zeta_1,\zeta_2) \mapsto
D^2_{\mathcal{G}}G_n(\varphi)(\zeta_1,\zeta_2)$ is continuous.
Proposition~\ref{prop:Ito} applied to the process $u_\lambda$ defined
by \eqref{strong} then yields
\begin{equation}
\label{ito}
\begin{split}
  &G_n(u_\lambda) + \int_0^\cdot \ip[\big]{Au_\lambda}{DG_n(u_\lambda)}(s)\,ds\\
  &\hspace{5em} = G_n(u_{0\lambda}) %
  + \int_0^\cdot DG_n(u_\lambda(s)) \bar{B}_\lambda(s)\,dW(s)\\
  &\hspace{5em} \quad + \int_0^\cdot \Bigl( DG_n(u_\lambda) \bar{F}_\lambda
    + \frac12 \operatorname{Tr}\bigl(\bar{B}_\lambda^*
    D^2_{\mathcal{G}}G_n(u_\lambda) \bar{B}_\lambda\bigr) \Bigr)(s)\,ds
\end{split}
\end{equation}
Recalling that $g'_n: \erre \to \erre$ is increasing, Lemma~\ref{lem:brez-str} implies that
\[
  \ip{Au_\lambda}{DG_n(u_\lambda)} = \ip{Au_\lambda}{g'(u_\lambda)}
  \geq 0,
\]
hence also, denoting a complete orthonormal system of $U$ by $(e_j)$,
\begin{align*}
  \int_\OO g_n(u_\lambda(t))
  &\leq \int_\OO g_n(u_{0\lambda})
    + \int_0^t g_n'(u_\lambda(s)) \bar{B}_\lambda(s)\,dW(s)\\
  &\quad + \int_0^t g_n'(u_\lambda(s)) \bar{F}_\lambda(s)\,ds
    + \frac12 \int_0^t \sum_{j=0}^\infty\int_\OO g_n''(u_\lambda(s))
    \abs[\big]{\bar{B}_\lambda(s)e_j}^2\,ds
\end{align*}
for every $t \in [0,T]$ and $n \in \enne$. We are now going to pass to
the limit as $n \to \infty$ in this inequality. Recalling that
$(g''_n)$ is uniformly bounded and that the paths of $u_\lambda$
belong to $C([0,T]; H)$ $\P$-a.s., the dominated convergence theorem yields
\begin{align*}
  \int_\OO g_n(u_\lambda(t))
  &\longrightarrow \frac12 \norm[\big]{u_\lambda^-(t)}^2
    \quad \forall t \in [0,T],\\
  \int_\OO g_n(u_{0\lambda})
  &\longrightarrow \frac12 \norm[\big]{u_{0\lambda}^-}^2.
\end{align*}
Note that $u_0$ is positive and $J_\lambda$ is positivity preserving,
hence $u_{0\lambda}=J_\lambda u_0$ is also positive and, in
particular, $u_{0\lambda}^-$ is equal to zero a.e. in $\OO$. Let us
introduce the (real) continuous local martingales
$(M^n)_{n \in \enne}$, $M$, defined as
\begin{align*}
  M^{\lambda,n}_t &:= \int_0^t g_n'(u_\lambda(s)) \bar{B}_\lambda(s)\,dW(s),\\
  M^\lambda_t &:= -\int_0^t (u_\lambda^-(s) \bar{B}_\lambda(s)\,dW(s).
\end{align*}
One has, by the ideal property of Hilbert-Schmidt operators,
\begin{align*}
  \bigl[M^{\lambda,n}-M^\lambda,M^{\lambda,n}-M^\lambda\bigr]_t
  &= \int_0^t \norm[\big]{(g'_n(u_\lambda(s)) + u_\lambda^-(s))%
    \bar{B}_\lambda(s)}_{\cL^2(U,\erre)}^2\,ds\\
  &\leq \int_0^t \norm[\big]{g'_n(u_\lambda(s)) + u_\lambda^-(s)}^2
    \norm[\big]{\bar{B}_\lambda(s)}^2_{\cL^2(U,H)}\,ds
\end{align*}
for all $t \in [0,T]$.  Recalling that
$u_\lambda\in L^0(\Omega;C([0,T];H))$ and $g_n'(r) \to -r^-$ for every
$r \in \erre$, it follows by the dominated convergence theorem that
$[M^{\lambda,n}-M^\lambda,M^{\lambda,n}-M^\lambda] \to 0$, hence that 
$M^{\lambda,n} \to M^\lambda$, as $n \to \infty$,
i.e. that
\[
  \int_0^t g_n'(u_\lambda(s)) \bar{B}_\lambda(s)\,dW(s)
  \longrightarrow -\int_0^t u_\lambda^-(s) \bar{B}_\lambda(s)\,dW(s)
\]
for all $t \in [0,T]$.
Similarly, the pathwise continuity of $u_\lambda$ and the dominated
convergence theorem yield
\[
  \int_0^t g_n'(u_\lambda(s)) \bar{F}_\lambda(s)\,ds
  \longrightarrow -\int_0^t u_\lambda^-(s) \bar{F}_\lambda(s)\,ds
\]
for all $t \in [0,T]$ as $n \to \infty$.
Finally, the pointwise convergence $g_n'' \to \ind{\erre_-}$ and the
dominated convergence theorem imply that
\[
\int_0^t \sum_{j=0}^\infty \int_\OO g_n''(u_\lambda(s))
\abs[\big]{\bar{B}_\lambda(s)e_j}^2\,ds
\longrightarrow \int_0^t \sum_{j=0}^\infty \int_\OO \ind{\{u_\lambda(s)<0\}}
\abs[\big]{\bar{B}_\lambda(s)e_j}^2\,ds
\]
for all $t \in [0,T]$ as $n \to \infty$.
We are thus left with
\begin{align*}
  \norm[\big]{u_\lambda^-(t)}^2
  &\leq \int_0^t \Bigl( -2\ip{u_\lambda^-(s)}{\bar{F}_\lambda(s)}
    + \sum_{j=0}^\infty \int_\OO \ind{\{u_\lambda(s)<0\}}
\abs[\big]{\bar{B}_\lambda(s)e_j}^2 \Bigr)\,ds\\
  &\quad -\int_0^t u_\lambda^-(s) \bar{B}_\lambda(s)\,dW(s).
\end{align*}
Let us now take the limit as $\lambda \to 0$: if follows from the
convergence property \eqref{conv} and the continuous mapping theorem that
\[
  \norm[\big]{u_\lambda^-(t)}^2 \longrightarrow \norm[\big]{\bar{u}^-(t)}^2.
\]
Recalling that $\bar{F}_\lambda = J_\lambda F(\bar{u})$, which
converges pointwise to $F(\bar{u})$, one has
\[
  \int_0^t -2\ip[\big]{u_\lambda^-(s)}{\bar{F}_\lambda(s)}\,ds
  \longrightarrow
  \int_0^t -2\ip[\big]{\bar{u}^-(s)}{F(s,\bar{u}(s))}\,ds
\]
Appealing again to \eqref{conv}, it is not difficult to check that
\[
  \ind{\{\bar{u}(s)<0\}} \leq \liminf_{\lambda\to0}
  \ind{\{u_\lambda(s)<0\}} \quad \text{a.e.~in } \OO \quad \forall s
  \in [0,T].
\]
Hence it follows from Fatou's lemma that
\[
\int_0^t \sum_{j=0}^\infty \int_\OO \ind{\{u_\lambda(s)<0\}}
\abs[\big]{\bar{B}_\lambda(s)e_j}^2\,ds \longrightarrow
\int_0^t \sum_{j=0}^\infty \int_\OO \ind{\{\bar{u}(s)<0\}}
\abs[\big]{B(s,\bar{u}(s))e_j}^2\,ds.
\]
Let us define the real continuous local martingales
$(M^\lambda)_{\lambda>0}$, $M$, defined as
\begin{align*}
M^\lambda_t &:= -\int_0^t u_\lambda^-(s) \bar{B}_\lambda(s)\,dW(s),\\
M_t &:= -\int_0^t \bar{u}^-(s) \bar{B}(s)\,dW(s).
\end{align*}
One has
\[
  \bigl[M^\lambda-M,M^\lambda-M\bigr]_t =
  \int_0^t \norm[\big]{u_\lambda^-(s)\bar{B}_\lambda(s)
    - \bar{u}^-(s)\bar{B}(s)}^2_{\cL^2(U,\erre)}\,ds,
\]
where, by the ideal property of Hilbert-Schmidt operators and the
contractivity of $J_\lambda$,
\begin{align*}
  \norm[\big]{u_\lambda^-\bar{B}_\lambda - \bar{u}^-\bar{B}}_{\cL^2(U,\erre)}
  &\leq \norm[\big]{(u_\lambda^- - \bar{u}^-)\bar{B}_\lambda}_{\cL^2(U,\erre)}
    + \norm[\big]{\bar{u}^-(\bar{B}_\lambda - \bar{B})}_{\cL^2(U,\erre)}\\
  &\leq \norm[\big]{u_\lambda^- - \bar{u}^-}
    \norm[\big]{\bar{B}}_{\cL^2(U,H)}
    + \norm[\big]{\bar{u}^-}
  \norm[\big]{\bar{B}_\lambda - \bar{B}}_{\cL^2(U,H)}.
\end{align*}
Since, as $\lambda \to 0$, $u_\lambda$ converges to $\bar{u}$ in the
sense of \eqref{conv} and, as already seen, $\bar{B}_\lambda$
converges to $\bar{B}$ in $L^0(\Omega;L^2(0,T;\cL^2(U,H)))$, the
dominated convergence theorem yields, for every $t \in [0,T]$,
\[
\bigl[M^\lambda-M,M^\lambda-M\bigr]_t \longrightarrow 0,
\]
thus also
\[
  \int_0^t u_\lambda^-(s) \bar{B}_\lambda(s)\,dW(s) \longrightarrow
  \int_0^t \bar u^-(s) B(s,\bar{u}(s))\,dW(s).
\]
Recalling assumption (A3), one obtains, for every $t \in [0,T]$,
\[ 
  \norm{\bar{u}^-(t)}^2 \leq 2C \int_0^t \norm{\bar{u}^-(s)}^2\,ds
  -2 \int_0^t \bar{u}^-(s) B(s,\bar{u}(s))\,dW(s),
\]
thus also, integrating by parts,
\[
  e^{-2Ct} \norm{\bar{u}^-(t)}^2 \leq 
  -2 \int_0^t e^{-2Cs} \bar{u}^-(s) B(s,\bar{u}(s))\,dW(s) =: \tilde M_t.
\]
The process $\tilde M$ is a positive local martingale, hence a
supermartingale, with $\tilde M(0)=0$, therefore $M$ is identically equal to
zero. This implies that $\norm{\bar{u}^-(t)}=0$ for all $t \in [0,T]$,
hence, in particular, that $u(t)$ is positive a.e. in $\OO$ for all
$t \in [0,T]$.
By definition of $\bar{u}$, we deduce that 
\[
  u^{\sigma} \geq 0 \quad\text{a.e.~in } \Omega \times[0,T] \times \OO
\]
for every $\sigma<\tau$. Since $\sigma$ is arbitrary, this readily
implies that
\[
  u \geq 0 \quad\text{a.e.~in } \co{0}{\tau} \times \OO,
\]
thus completing the proof of Theorem~\ref{thm:pos}.


\section{Positivity of forward rates}
Musiela's stochastic PDE can be written as
\begin{equation}
  \label{eq:Mus}
  du + Au\,dt = \beta(t,u)\,dt
  + \sum_{k=1}^\infty \sigma_k(t,u)\,dw^k(t), \qquad u(0)=u_0,
\end{equation}
where $-A$ is (formally, for the moment) the infinitesimal generator
of the semigroup of translations, $(w^k)_{k\in\enne}$ is a sequence of
independent standard Wiener processes, $\sigma_k$ is a random,
time-dependent superposition operator for each $k \in \enne$, as well
as $\beta$, and $u$ takes values in a space of continuous functions,
so that $u(t,x):=[u(t)](x)$, $x \geq 0$, models the value of the
forward rate prevailing at time $t$ for delivery at time $t+x$. In
order to exclude arbitrage (or, more precisely, in order for the
corresponding discounted bond price process to be a local martingale),
$\beta$ needs to satisfy the so-called Heath-Jarrow-Morton
no-arbitrage condition
\[
  \beta(t,v) = \sum_{k=1}^\infty \sigma_k(t,v)
  \int_0^\cdot [\sigma_k(t,v)](y)\,dy.
\]
In order for \eqref{eq:Mus} to admit a solution with continuous paths,
a by now standard choice of state space is the Hilbert space
$H_\alpha$, $\alpha>0$, which consists of absolutely continuous
functions $\phi:\erre_+ \to \erre$ such that
\[
  \norm[\big]{\phi}^2_{H_\alpha} := \phi(\infty)^2
  + \int_0^\infty \abs{\phi'(x)}^2 e^{\alpha x}\,dx < \infty.
\]
Under measurability, local boundedness, and local Lipschitz continuity
conditions on $(\sigma_k)$, one can rewrite \eqref{eq:Mus} as
\begin{equation}
  \label{eq:Musa}
  du + Au\,dt = \beta(t,u)\,dt + B(t,u)\,dW(t), \qquad u(0)=u_0,  
\end{equation}
where $A$ is the generator of the semigroup of translations on
$H_\alpha$, $W$ is a cylindrical Wiener process on $U=\ell^2$, and
$B:\Omega \times \erre_+ \times H \to \cL^2(U,H)$ is such that
\[
  \sum_{k=1}^\infty \int_0^\cdot \sigma_k(s,v(s))\,dw^k(s) =
  \int_0^\cdot B(s,v(s))\,dW(s).
\]
Under such assumptions on $(\sigma_k)$, \eqref{eq:Mus} admits a unique
local mild solution with values in $H_\alpha$. If $(\sigma_k)$ satisfy
stronger (global) boundedness and Lipschitz continuity assumptions,
then the local mild solution is in fact global.  For details we refer
to \cite{filipo}, as well as to \cite{cm:pos1}.

Positivity of forward rates, i.e. of the mild solution to
\eqref{eq:Mus}, is established in \cite{cm:pos1} by proving positivity
of mild solutions in weighted $L^2$ spaces to regularized versions of
\eqref{eq:Mus}. Such an approximation argument is employed because the
conditions on $(\sigma_k)$ ensuring (local) Lipschitz continuity of
the coefficients in the associated stochastic evolution
\eqref{eq:Musa} equation in $H_\alpha$ do not imply (local) Lipschitz
continuity of the coefficients if state space is changed to a weighted
$L^2$ space.

Thanks to Theorem~\ref{thm:pos}, we can give a much shorter, more
direct proof of the (criterion for the) positivity of forward
rates. Let $L^2_{-\alpha}$ denote the weighted space
$L^2(\erre_+,e^{-\alpha x}\,dx)$, and note that $H_\alpha$ is
continuously embedded in $L^2_{-\alpha}=:H$.
Let us check that assumptions (A1), (A2), and (A3) are satisfied.
Assumption (A1) holds true with the choice $\OO=\erre_+$, endowed with
the absolutely continuous measure $m(dx):=e^{-\alpha x}\,dx$.  As far
as assumption (A2) is concerned, a simple computation shows that
$A+\alpha I$ is monotone on $L^2_{-\alpha}$, and, by standard ODE
theory, one also verifies that the range of $A + \alpha I + I$
coincides with the whole space $L^2_{-\alpha}$, therefore $A+\alpha I$
is maximal monotone. Even though $A$ itself is not maximal monotone,
this is clearly not restrictive, as the ``correction'' term $\alpha I$
can be incorporated in $\beta$ without loss of generality. To verify
that the resolvent $J_\lambda \in \cL(H)$ of $A+\alpha I$ is
sub-Markovian, let $y \in H$, so that $J_\lambda y \in \dom(A)$ is the
unique solution $y_\lambda$ to the problem
\[
  y_\lambda - \lambda y_\lambda' + \lambda\alpha y_\lambda = y.
\]
If $0\leq y\leq1$ a.e. in $\erre_+$, then we have, multiplying both
sides by ${(y_\lambda-1)^+}$, in the sense of the scalar product of $H$, that
\begin{align*}
  &(1+\lambda\alpha) \ip[\big]{y_\lambda}{(y_\lambda-1)^+}_{-\alpha}
  - \lambda \ip[\big]{y'_\lambda}{(y_\lambda-1)^+}_{-\alpha}\\
  &\hspace{5em}= \ip[\big]{y}{(y_\lambda-1)^+}_{-\alpha}
    \leq \ip[\big]{1}{(y_\lambda-1)^+}_{-\alpha}.
\end{align*}
Here and in the following we denote the scalar product and norm of
$L^2_{-\alpha}$ simply by $\ip{\cdot}{\cdot}_{-\alpha}$ and
$\norm{\cdot}_{-\alpha}$, respectively.
Since
\begin{equation}
  \label{eq:iole}
  \ip[\big]{y_\lambda}{(y_\lambda-1)^+}_{-\alpha} -
  \ip[\big]{1}{(y_\lambda-1)^+}_{-\alpha} =
  \norm[\big]{(y_\lambda-1)^+}^2_{-\alpha},
\end{equation}
we obtain
\[
  \norm[\big]{(y_\lambda-1)^+}^2_{-\alpha}
  - \frac{\lambda}{2} \int_0^{\infty} \frac{d}{dx}((y_\lambda-1)^+)^2(x)
  e^{-\alpha x}\,dx
  + \lambda\alpha \ip[\big]{y_\lambda}{(y_\lambda-1)^+}_{-\alpha}
  \leq 0,
\]
where, integrating by parts,
\begin{align*}
  &- \frac{\lambda}{2} \int_0^{\infty} \frac{d}{dx}((y_\lambda-1)^+)^2(x)
    e^{-\alpha x}\,dx\\
  &\hspace{5em} = -\frac{\lambda\alpha}{2} \int_0^{\infty}
    ((y_\lambda(x)-1)^+)^2 e^{-\alpha x}\,dx
    + \frac{\lambda}{2} ((y_\lambda(0)-1)^+)^2\\
  &\hspace{5em} = -\frac{\lambda\alpha}{2}
    \ip[\big]{y_\lambda}{(y_\lambda-1)^+}_{-\alpha}
    + \frac{\lambda\alpha}{2} \ip[\big]{1}{(y_\lambda-1)^+}_{-\alpha}
    + \frac{\lambda}{2} ((y_\lambda(0)-1)^+)^2\\
  &\hspace{5em} \geq -\frac{\lambda\alpha}{2}
    \ip[\big]{y_\lambda}{(y_\lambda-1)^+}_{-\alpha}.
\end{align*}
Rearranging the terms yields
\[
  \norm[\big]{(y_\lambda-1)^+}^2_{-\alpha}
  + \frac{\lambda\alpha}{2} \ip[\big]{y_\lambda}{(y_\lambda-1)^+}_{-\alpha}
  \leq 0,
\]
where the second term on the left-hand side is positive by
\eqref{eq:iole}. Therefore $\norm{(y_\lambda-1)^+}_{-\alpha}=0$, which
implies that $y_\lambda \leq 1$ a.e. in $\erre_+$. A completely
similar argument, i.e. scalarly multiplying the resolvent equation by
$y_\lambda^-$, also shows that $y_\lambda \geq 0$ a.e. in $\erre_+$,
thus completing the proof that $J_\lambda$ is sub-Markovian.
We still need to show that $J_\lambda$ is contractive in
$L^1_{-\alpha}$.  Let $y,z \in H$ and $y_\lambda:=J_\lambda y$,
$z_\lambda:=J_\lambda z$, so that
\begin{equation}
  \label{eq:yzl}
  (y_\lambda-z_\lambda) - \lambda(y_\lambda - z_\lambda)'
  + \lambda\alpha(y_\lambda-z_\lambda) 
  = y-z.
\end{equation}
Define the sequences of functions
$(\gamma_k), (\hat{\gamma}_k) \subset \erre^\erre$ as
\[
  \gamma_k: r \mapsto \tanh(kr), \qquad
  \hat{\gamma}_k: r \mapsto \int_0^r \gamma_k(s)\,ds,
\]
and recall that, as $k \to \infty$, $\gamma_k$ converges pointwise to
the sign function, and $\hat{\gamma}_k$ converges pointwise to the
absolute value function. Scalarly multiplying \eqref{eq:yzl} with
$\gamma_k(y_\lambda-z_\lambda)$ yields
\begin{align*}
  &(1+\lambda\alpha)
  \ip[\big]{y_\lambda-z_\lambda}{\gamma_k(y_\lambda-z_\lambda)}_{-\alpha}
  - \lambda \ip[\big]{(y_\lambda-z_\lambda)'}%
    {\gamma_k(y_\lambda-z_\lambda)}_{-\alpha}\\
  &\hspace{5em} = \ip[\big]{y-z}{\gamma_k(y_\lambda-z_\lambda)}_{-\alpha}
    \leq \norm[\big]{y-z}_{L^1_{-\alpha}},
\end{align*}
where, integrating by parts,
\begin{align*}
  \ip[\big]{(y_\lambda-z_\lambda)'}{\gamma_k(y_\lambda-z_\lambda)}_{-\alpha}
  &= \int_0^\infty \bigl( \gamma_k(y_\lambda-z_\lambda)(x)
    (y_\lambda-z_\lambda)'(x) \bigr) e^{-\alpha x}\,dx\\
  &= \int_0^\infty \frac{d}{dx} \hat{\gamma}_k(y_\lambda-z_\lambda)(x)
    e^{-\alpha x}\,dx\\
  &= - \hat{\gamma}_k(y_\lambda(0)-z_\lambda(0))
    + \alpha \int_0^\infty \hat{\gamma}_k(y_\lambda-z_\lambda)(x)
    e^{-\alpha x}\,dx\\
  &\leq \alpha \int_0^\infty \hat{\gamma}_k(y_\lambda-z_\lambda)(x)
    e^{-\alpha x}\,dx.
\end{align*}
This implies
\begin{align*}
  &\ip[\big]{y_\lambda-z_\lambda}{\gamma_k(y_\lambda-z_\lambda)}_{-\alpha}\\
  &\hspace{3em} + \lambda\alpha
    \ip[\big]{y_\lambda-z_\lambda}{\gamma_k(y_\lambda-z_\lambda)}_{-\alpha}
    - \lambda\alpha \int_0^\infty \hat{\gamma}_k(y_\lambda-z_\lambda)(x)
    e^{-\alpha x}\,dx\\
  &\hspace{5em} \leq \norm[\big]{y-z}_{L^1_{-\alpha}}.
\end{align*}
Taking the limit as $k \to \infty$, the sum of the second and third
term on the left-hand side converges to zero by the dominated
convergence theorem, while the first term on the left-hand side
converges to $\norm{y_\lambda-z_\lambda}_{L^1_{-\alpha}}$, thus
proving that
\[
  \norm[\big]{y_\lambda-z_\lambda}_{L^1_{-\alpha}} \leq
  \norm[\big]{y-z}_{L^1_{-\alpha}},
\]
i.e. that the resolvent of $A+\alpha I$ is contractive in
$L^1_{-\alpha}$. We have thus shown that assumption (A2) holds for
$A+\alpha I$.
Moreover, assumption (A3) is satisfied if, for example, 
\[
\abs{\sigma_k(\omega,t,x,r)} \ind{\{r \leq 0\}} \lesssim r^-.
\]
for all $k \in \enne$ and $(\omega,t,x) \in \Omega \times \erre_+^2$
(see \cite{cm:pos1}, where also slightly more general sufficient
conditions are provided).
Since all integrability assumptions of Theorem~\ref{thm:pos} are
satisfied, as it follows by inspection of the proof of well-posedness
in $H_\alpha$ (see~\cite{filipo,cm:pos1,cm:MF10}), we conclude that,
under the above assumptions on $(\sigma_k)$, forward rates are
positive at all times.

\bibliographystyle{amsplain}
\bibliography{ref,extra}

\end{document}